\documentstyle[12pt]{article}
\newtheorem{theor}{Theorem}

\newtheorem{lem}{Lemma}
\begin{document} 
\title{Rank 1 forms, closed zones and laminae}
\author{Michel Deza\\
Ecole Normale Sup\'erieure, Paris \\ \and 
Viatcheslav Grishukhin\\
CEMI, Russian Academy of Sciences, Moscow}
\date{} 
\maketitle 

\begin{abstract}
For a given lattice, we establish an equivalence involving a closed 
zone of the corresponding Voronoi polytope, a lamina hyperplane of 
the corresponding Delaunay partition and a rank 1 quadratic form
being an extreme ray of the corresponding $L$-type domain.
\end{abstract} 

An $n$-dimensional lattice determines two normal partitions of 
the $n$ space ${\bf R}^n$ into polytopes. 
These are the Voronoi partition and the Delaunay partition. These 
partitions are dual, i.e. a $k$-dimensional face
of one partition is orthogonal to an $(n-k)$-dimensional face of 
the other partition. Besides, a vertex of one partition is the 
center of a polytope of the other partition. 

The Voronoi partition consists of Voronoi polytopes with its centers 
in lattice points. Moreover, any polytope of the Voronoi partition
is obtained by a translation of the Voronoi polytope with the center 
in the origin (=the zero lattice point). Call this polytope {\em the}
Voronoi polytope. It consists of those points of ${\bf R}^n$ that are
at least as closed to 0 as to any other lattice point.

The Delaunay partition consists of Delaunay polytopes which are, in
general, not congruent. The set of all Delaunay polytopes having 
0 as a vertex is called the {\em star} of Delaunay polytopes.
Each Delaunay polytope is the convex hull of all lattice points
lying on an {\em empty} sphere.  This sphere is called empty, since
no lattice point is an interior point of the sphere.

The Voronoi polytope and the Delaunay polytopes of the star are 
tightly related to minimal vectors of cosets $2L$ in $L$. A coset $Q$ 
is called {\em simple} if it contains, up to sign, only one minimal vector. 

For a Delaunay polytope $P_D$, the lattice vector between any two 
vertices of $P_D$ is a minimal vector of a coset of $L/2L$. 
A lattice vector is an edge of a Delaunay polytope of the star (and 
then, by duality, it defines a facet of the Voronoi polytope) if and 
only if it is the minimal vector of a simple coset of $L/2L$. All 
minimal vectors of a non-simple coset are diagonals of a symmetric face 
of a Delaunay polytope of the star. 

The set ${\cal P}(P)$ of all faces of all dimensions of a polytope $P$ 
is  partially ordered by inclusion. Call it {\em face poset of $P$}. 
The face poset of the Voronoi polytope $P_V$ of a lattice $L$ determines
uniquely the combinatorial structure of $P_V$ and the {\em $L$-type}
of the Voronoi and Delaunay partitions. The notion of an $L$-type was
introduced by G.Voronoi in \cite{Vo}. One says that a lattice $L$
{\em belongs to} or {\em is of} an $L$-type if its Voronoi partition
has this $L$-type. In other words, two lattices (and their Voronoi
and Delaunay partitions) belong to the same $L$-type if the
corresponding partitions are combinatorially and topologically
equivalent, or, equivalently, the face posets of their Voronoi
polytopes are isomorphic.
If we reverse the order of ${\cal P}(P_V)$, we obtain the poset of
those faces of Delaunay polytopes of the star that contain the point 0.

There is a small perturbation of a basis of $L$ that does not change
the $L$-type of $L$. Suppose that a perturbation of the basis changes
the $L$-type. Then the Delaunay partition changes, and there is an
empty sphere such that a lattice point either leaves or comes onto the
sphere.

There is the following simple but important test of emptiness of a
sphere. Let $S \subset {\bf R}^n$ be an $(n-1)$-dimensional sphere
with the origin point 0 on it. Let $v_1,v_2,...,v_n$ be $n$ linearly
independent lattice vectors with endpoints on $S$. Let
$u \in {\bf R}^n$ be an arbitrary vector and $u=\sum_{i=1}^nz_iv_i$.
(We denote by $(p,q)$ the scalar product of vectors $p$ and $q$, and
set $p^2=(p,p)$).
\begin{lem}
\label{uS}
{\rm (Proposition 4 of \cite{BG}).} The endpoint of $u$ is not an
interior point of $S$ if and only if the following inequality holds
\begin{equation}
\label{uv}
u^2 \ge \sum_{i=1}^n z_iv_i^2.
\end{equation}
The endpoint of $u$ lies on $S$ if and only if (\ref{uv}) holds as
equality.
\end{lem}
{\bf Proof}. Let $c$ be the center of $S$. Since the endpoint of
$v_i$ lies on $S$, we have $(v_i-c)^2=c^2$, $1 \le i \le n$, i.e.
$v_i^2=2(v_i,c)$. Multiplying this equality by $z_i$, summing over $i$
and taking in attention that $u=\sum z_iv_i$, we obtain
\[\sum_{i=1}^nz_iv_i^2=2(u,c). \]
Since the endpoint of $u$ is not an interior point of $S$,
$(u-c)^2\ge c^2$, i.e. $u^2\ge 2(u,c)$. Using the above equality, we
obtain (\ref{uv}). It is easy to see that (\ref{uv}) holds as
equality if and only if the inequality $u^2\ge 2(u,c)$ holds as
equality, i.e. if and only if the endpoint of $u$ lies on $S$.
\hfill $\Box$

Any basis ${\cal B}=\{b_i,1\le i\le n\}$ of an $n$-dimensional
lattice $L$ determines uniquely a positive definite quadratic form
\[f(x)=(\sum_1^n b_ix_i)^2=\sum_{1 \le i,j \le n}a_{ij}x_ix_j. \]
The symmetric matrix $a_{ij}$ of the coefficients of this form
is the Gram matrix of the basis $\cal B$, i.e. $a_{ij}=(b_i,b_j)$.
The matrix $a_{ij}$ can be considered as a point of
an $N$-dimensional space, where $N={n+1 \choose 2}$. In this space, all
positive definite forms form an open cone. The closure of this cone is
the cone ${\cal P}_n$ of all positive semi-definite quadratic forms of
order $n$.

One says that a quadratic form belongs to or is of an $L$-type if its 
lattice belongs to this $L$-type. Hence the cone ${\cal P}_n$ is 
partitioned into $L$-type domains of forms of the same $L$-type. Of 
course, the cone ${\cal P}_n$ has many domains of the same $L$-type 
corresponding to distinct choices of a basis. Voronoi proved that each
$L$-type domain is an open polyhedral cone of dimension $k$, 
$1 \le k \le N$. An $N$-dimensional $L$-type domain is called
{\em general}. Domains of other dimensions are called {\em special}.
Any face of the closure of a general $L$-type domain is the closure of
a special $L$-type domain. One-dimensional $L$-type domains are extreme
rays of the closure of a general $L$-type domain.

Call a form an {\em edge form} if it belongs to a one-dimensional 
$L$-type domain. In \cite{BG}, an edge form is called {\em rigid form},
since the only transformation of the corresponding lattice that
does not change its $L$-type is a scaling.
A typical edge form is the square of a linear form:
$f(x)=(\sum_{i=1}^n p_ix_i)^2$, i.e. a rank 1 form. But, for $n \ge 4$,
there are edge forms of full rank $n$.

A polyhedral domain of quadratic forms is called {\em dicing domain} if
all extreme rays of its closure are forms of rank 1. Dicings were defined
and studied by Erdahl and Ryshkov \cite{ER}. In \cite{ER}, they give
conditions when a dicing domain is an $L$-type domain and prove the
following theorem (Theorem 4.3 of \cite{ER}): {\em An $L$-type domain
is a dicing domain if and only if all the edge forms are rank 1 forms}.

Return to the Voronoi polytope $P_V$ of an $n$-dimensional lattice $L$.
The Voronoi polytope $P_V$ itself and its facets (i.e. faces of dimension
$n-1$) are centrally symmetric, so as its vertices and edges. But 
faces of other dimensions are, in general, not centrally symmetric. 
For example, there are many types of two-dimensional faces: hexagons, 
and others that are degenerated cases of a hexagon when some of its 
edges are compressed to a point. 

The set of edges of $P_V$ is partitioned into classes of mutually parallel
edges. These classes are called {\em zones}. There are two types of zones:
closed and open. A zone is called {\em closed} if every two-dimensional 
face contains either two edges of the zone or else none. Otherwise the
zone is called {\em open}. The notions of closed and open zones was 
introduced by P.Engel in \cite{En}. 

A closed zone has the following property. Let $l$ be the minimal length 
of edges of a closed zone $Z$. Let us shorten all edges of $Z$ onto a 
value $\varepsilon \le l$. If $\varepsilon<l$, then $Z$ 
remains closed, and the new polytope $P'_V$ (with shortened edges)
is a Voronoi polytope with the same face poset as $P_V$. If $\varepsilon=
l$, then $Z$ transforms into an open zone, and $P'_V$ has another face
poset, since at least one edge vanishes.

Since the Voronoi partition is dual to the Delaunay partition, each 
edge of a Voronoi polytope is orthogonal to a facet of a Delaunay 
polytope. A facet $F$ of a Delaunay polytope of a lattice $L$ generates
an affine hyperplane $H$ in ${\bf R}^n$, namely the hyperplane, where
$F$ lies. Obviously $F$ contains $n$ affinely independent lattice points.
Hence the intersection $L\cap H$ is an $(n-1)$-dimensional sub-lattice of
$L$. The Delaunay partition of $L$ generates a partition of the hyperplane
$H$ into Delaunay polytopes of the lattice $L\cap H$. It may be that all
$(n-1)$-dimensional Delaunay polytopes of the partition of $H$ are facets
of Delaunay polytopes of the original Delaunay partition of $L$.
In this case $H$ is called a {\em lamina} of the lattice $L$.

The notion of lamina was introduced and extensively used by Ryshkov and
Baranovskii in \cite{RB} (see \S9.4). If $L$ belongs to a general
$L$-type and if a hyperplane $H$ is not a lamina, then it intersects
in an interior point an edge of at least one Delaunay polytope $P_D$ of
the star (Lemma 9.3 of \cite{RB}). In other words, there are two vertices
of $P_D$ that lie in distinct half-spaces determined by $H$. We reformulate
Lemma 9.3 of \cite{RB} for a lattice of an arbitrary $L$-type.
\begin{lem}
\label{lRB}
A hyperplane is a lamina if and only if it does not intersect any
Delaunay polytope of the star in an interior point.
\end{lem}

Obviously, a lamina determines a family of parallel laminae that
partitions the lattice $L$ into parallel layers, each of them spanning
a lamina. Lemma~\ref{lRB} implies the following corollary.

{\bf Corollary} {\em Every Delaunay polytope of a lattice with a lamina
lies between two neighboring laminae with vertices on these two laminae.}

The main property of this partition of $L$ into layers, spanning laminae,
is that the distances between layers may be changed without changing
the $L$-type of $L$.

Let $H$ be a hyperplane spanning an $(n-1)$-dimensional sub-lattice of $L$.
Let $e$ be a unit vector orthogonal to $H$. Then we can define an
$\epsilon$-extension along $e$ of the space and of the lattice $L$ as
follows. Any vector $v \in {\bf R}^n$ is uniquely decomposed as
$v=v_e+v_H$, where $v_e=(e,v)e$ and $v_H$ are the projections of $v$
onto $e$ and $H$, respectively. An $\epsilon$-extension of ${\bf R}^n$
along $e$ transforms every vector $v$ into the vector
\begin{equation}
\label{ext}
v'=(1+\epsilon)v_e+v_H=\epsilon(e,v)e+v.
\end{equation}
Here the $\epsilon$-extension is in fact a contraction if $\epsilon<0$.
In particular, for the norm (=squared length) ${v'}^2$ of the extended
vector, we obtain the following expression, where we set
$\lambda=\epsilon(2+\epsilon)$:
\begin{equation}
\label{v2}
{v'}^2=v^2+\lambda(e,v)^2.
\end{equation}
Of course, an $\epsilon$-extension of ${\bf R}^n$ along $e$ transforms
a lattice $L \subset {\bf R}^n$ into an {\em extended lattice}
$L^{\epsilon}$.

There is the following relation between the above introduced notions 
of a rank 1 form, a closed zone, a lamina and the extended lattice.

\begin{theor}
Let $L$ be an $n$-dimensional lattice, $H$ be a hyperplane spanning
an $(n-1)$-dimensional sub-lattice of $L$, $e$ be an $n$-dimensional
unit vector orthogonal to $H$. Let $f(x)$ be the quadratic form
corresponding to a basis $\{b_i:1 \le i \le n\}$ of $L$ and $D(f)$ be
the $L$-type domain of $f$. The following assertions are equivalent:

(i) $H$ is a lamina of the Delaunay partition of $L$;

(ii) the Voronoi polytope $P_V$ of $L$ has a closed zone $Z_e$ of edges
parallel to the vector $e$;

(iii) the $\epsilon$-extended along $e$ lattice $L^{\epsilon}$ has the
same $L$-type as $L$ for all $\epsilon>0$;
 
(iv) the rank 1 form $f_e(x)=(e,\sum_1^n b_ix_i)^2$ lies on an extreme
ray of the closure of $D(f)$, i.e. $f+\lambda f_e \in D(f)$ for all
nonnegative $\lambda$.
\end{theor}
{\bf Proof}. (i)$\Rightarrow$(ii). Let $H$ be a lamina partitioned into 
facets of Delaunay polytopes of $L$. We can suppose that $H$ contains 
the origin 0. Consider the edges of the Voronoi polytope $P_V$ that are
orthogonal to the lying in $H$ facets of the star. Obviously these edges
are parallel to $e$ and form the zone $Z_e$. We show that $Z_e$ is closed.
If not, there is a 2-face $T$ of $P_V$ containing exactly one edge
$u_1 \in Z_e$. The edges of $T$ form a polygon. Let $u_1, u_2,...,u_k$
be consecutive edges of this polygon. Let $F_i$ be the facet of the star
that is orthogonal to the edge $u_i$, $1 \le i \le k$. The set of facets
$\{F_i:1 \le i \le k\}$ has a common $(n-2)$-dimensional face of the star
that is orthogonal to $T$ and lies in the lamina $H$. Since, for
$2\le i \le k$, the edge $u_i$ is not parallel to $u_1$, the facet $F_i$
does not lie in the lamina $H$. Hence there is an index $j$ such that
the facets $F_j$ and $F_{j+1}$ lie in distinct halfspaces separated by
$H$. Obviously, $F_j$ and $F_{j+1}$ are facets of a same Delaunay polytope
$P_j$ of the star, and $H$ intersects $P_j$.
This contradicts to definition of a lamina. Hence $Z_e$ cannot be open.

(ii)$\Rightarrow$(i). Let $u_1 \in Z_e$. Consider the facet $F_1$ of
the star that is orthogonal to $u_1$ and contains $0 \in L$. $F_1$ spans 
a hyperplane $H$ that is orthogonal to $e$ and contains $0 \in L$.
Let $T_1$ be a 2-face of $P_V$ containing $u_1$, and $u_2$ be the
second edge from $Z_e$ contained in $T_1$. Let $F_2$ be the facet of
the star that is dual to $u_2$. $F_2$ intersects $F_1$ by an $(n-2)$-face
that is dual to $T_1$ and contains 0. Hence $F_2$ contains 0 and is 
orthogonal to $e$. This implies that $F_2$ lies in $H$. Similarly,
for $i=3,4,...$, we consider the 2-face $T_i$ containing 
$u_i,u_{i+1} \in Z_p$, and prove that the facet $F_i$ of the star 
dual to $u_i$ lies in $H$. Since $P_V$ is a polytope, it has a finite
number of edges. Hence there is $i_0$ such that $u_{i_0}=u_1$. 
We obtain a set of facets of the star lying in $H$ such that the facet 
$F_i$ intersects $F_{i-1}$ and $F_{i+1}$. Therefore the
intersection of the star with $H$ consists of facets of the star. 
Since this is true for all Voronoi polytopes having centers in $H$, 
the hyperplane $H$ is partitioned into facets of Delaunay polytopes. 
This means that $H$ is lamina.  

(i)$\Rightarrow$(iii). 
It is sufficient to prove that no lattice point comes onto or leaves
the empty sphere $S$ of a Delaunay polytope $P_D$. Without loss of
generality, we can suppose that $P_D$ belongs to the star. Let
$v_1, v_2,...,v_n$ be $n$ linearly independent lattice vectors with
endpoint in vertices of $P_D$, i.e. they lie on the sphere $S$.
Let $u$ be any lattice vector and let $u=\sum_{i=1}^nz_iv_i$ be its
decomposition by $v_i$, $1 \le i \le n$. According to (\ref{ext}),
after an $\epsilon$-extension of the space along $e$, $u$ and $v_i$
are transformed into the vectors
\[u'=u+\epsilon (u,e)e, \mbox{  }v'_i=v_i+\epsilon (v_i,e)e. \]
By Lemma~\ref{uS}, it is sufficient to prove that the inequality
${u'}^2 \ge \sum_{i=1}^n z_i{v'_i}^2$ is strict or is an equality
according to the inequality $u^2 \ge \sum_{i=1}^n z_iv_i^2$ is strict 
or is an equality.

Recall that $H$ contains an ($n-1$)-dimensional sub-lattice. Hence
$(v,e)=k(v)\alpha$ for any lattice vector $v$, where $k(v)$ is
an integer and $\alpha$ does not depend on $v$. Let $k=k(u)$ and
$k_i=k(v_i)$. Multiplying the equality $u=\sum_{i=1}^nz_iv_i$ by $e$,
we obtain the equality $k=\sum_{i=1}^nz_ik_i$. By Corollary of
Lemma~\ref{lRB}, the vertices of $P_D$ lie on two neighboring laminae.
Hence $|k_i|=0,1$. Without loss of generality we can suppose that
$k_i\ge 0$, i.e. $k_i=0,1$ and therefore $(v_i,e)^2=k_i^2=k_i$. Hence,
using (\ref{v2}), we have
\[{v'_i}^2=v_i^2+\lambda (v_i,e)^2=v_i^2+\lambda k_i \alpha^2. \]
Since $\sum_{i=1}^n z_ik_i=k$, we obtain the equality
\begin{equation}
\label{uzv}
\sum_{i=1}^nz_i{v'_i}^2=\sum_{i=1}^n z_iv_i^2+
\lambda \alpha^2 \sum_{i=1}^n z_ik_i=
\sum_{i=1}^n z_i v_i^2+\lambda\alpha^2 k.
\end{equation}
Suppose that $u^2>\sum_{i=1}^n z_iv_i^2$. We show that then
$\sum_1^n z_i{v'_i}^2<{u'}^2=u^2+\lambda(u,e)^2=u^2+\lambda \alpha^2 k^2$.
For an integer $k$, we have $k \le k^2$ with equality if and only if
$k=0,1$. The above inequalities and the equality (\ref{uzv}) imply
\[\sum_1^n z_i{v'_i}^2<u^2+\lambda \alpha^2 k^2={u'}^2. \]
Now let $u^2=\sum_1^n z_iv_i^2$. This means that $u$ has the endpoint
in a vertex of $P_D$. Hence $k=0,1$, i.e. $k^2=k$. Hence the equality
(\ref{uzv}) implies the equality
\[ \sum_1^n z_i{v'_i}^2=u^2+\lambda \alpha^2 k^2={u'}^2. \]

(iii)$\Rightarrow$(iv).
We prove in the implication (i)$\Rightarrow$(iii) that
$f(L^{\epsilon}) \in D(f)$ for every $\epsilon \ge 0$. Now we show that
\[f(L^{\epsilon})=f(L)+\lambda f_e, \]
where $\lambda=\epsilon(2+\epsilon)$.
The basic vectors of the extended lattice $L^{\epsilon}$ have the form
\[b'_i=b_i+\epsilon(e,b_i)e, \hspace{3mm}1 \le i \le n. \]
Hence the coefficients $a'_{ij}$ of the quadratic form
$f^{\epsilon}=f(L^{\epsilon})$ are as follows
\[a'_{ij}=(b'_i,b'_j)=(b_i,b_j)+\lambda(e,b_i)(e,b_j). \]
Hence we obtain
\[f^{\epsilon}(x)=f(x)+\lambda(\sum_{i=1}^n(e,b_i)x_i)^2=
f(x)+\lambda f_e(x). \]
This means that the ray $\{\lambda f_e: \lambda \ge 0\}$ belongs to
the closure of $D(f)$. Since this ray is a one-dimensional $L$-type
domain, it is an extreme ray of cl$D(f)$.

(iv)$\Rightarrow$(iii).
If $f_e$ lies on an extreme ray of cl$D(f)$, then the quadratic 
function $f^{\epsilon}=f+\epsilon(2+\epsilon)f_e$ belongs to $D(f)$ 
for all $\epsilon \ge 0$. The matrix $a'_{ij}$ of $f^{\epsilon}$ is
\[a_{ij}+\epsilon(2+\epsilon)(e,b_i)(e,b_j)=(b'_i,b'_j). \]
Hence $f^{\epsilon}$ is a quadratic form of the $\epsilon$-extended
lattice $L^{\epsilon}$. So, $L^{\epsilon}$ has the same $L$-type as
$L$ for all $\epsilon \ge 0$.

(iii)$\Rightarrow$(i).
We show that if the extended along $e$ lattice $L^{\epsilon}$ has the
same $L$-type as $L$ for all $\epsilon>0$, then the hyperplane $H$
which is orthogonal to $e$ is a lamina. Suppose $H$ is not a lamina.
Then $H$ intersects a Delaunay polytope $P_D$ of the star in an interior
point. Hence there are two vectors $v_1$ and $v_2$ with endpoints in
vertices of $P_D$ such that $k_1=k(v_1)\ge 1$ and $k_2=k(v_2) \le -1$.
Consider the lattice vector
\[u=q(k_1v_2-k_2v_1), \]
where the integer $q$ is chosen such that the endpoint of $u$ does not
lie on the empty sphere $S$ circumscribing $P_D$. Expand the pair of
vectors $v_1$, $v_2$ up to a set of $n$ independent vectors with
endpoints in vertices of $P_D$. Then the above expression for $u$ is
the representation of $u$ as a linear combination of these $n$ 
independent vectors. Since the endpoint of $u$ does not lie on $S$, by
Lemma~\ref{uS}, $\Delta \equiv u^2-q(k_1v_2^2-k_2v_1^2)>0$.

Consider an $\epsilon$-extension of the space along $e$. Let $u'$, $v'_1$,
$v'_2$ be the $\epsilon$-extended vectors. Consider the difference
$\Delta'={u'}^2-q(k_1{v'_2}^2-k_2{v'_1}^2)$. Using (\ref{v2}), we
obtain
\[\Delta'=\Delta+\lambda(e,u)^2-q(k_1\lambda(e,v_2)^2-
k_2\lambda(e,v_1)^2). \]
Since $(e,u)=\alpha k(u)=0$, $(e,v_1)=\alpha k_1$,
$(e,v_2)=\alpha k_2$ and $k_2<0$, we have
\[\Delta'=\Delta-\lambda qk_1 |k_2|(|k_2|+k_1). \]
Let $p=qk_1 |k_2|(k_1+|k_2|)>0$. Then for $\lambda=\frac{\Delta}{p}$,
we obtain $\Delta'=0$. Lemma~\ref{uS} implies that, for this $\lambda$,
the endpoint of $u$ lies on $S$. This means that there is $\epsilon>0$
such that the $\epsilon$-extended lattice $L^{\epsilon}$
has the $L$-type distinct from the $L$-type of $L$. We obtain a
contradiction.
\hfill $\Box$

\vspace{3mm}
{\bf Remark 1}.
The rank 1 function $f_e(x)=(e,\sum_1^n b_ix_i)^2$ depends on the
basis corresponding to the function $f(x)$. But this dependence is
not essential in the following sense. Recall that
$(e,b_i)=\alpha k(b_i)$, where $k_i=k(b_i)$ is an integer. If we slightly
move $f(x)$ in the domain $D(f)$, then we slightly change the basis
${\cal B}=\{b_i:1 \le i \le n\}$ and the scalar products $(e,b_i)=
\alpha k_i$. But, since $k_i$ is an integer, it cannot change slightly.
Hence only $\alpha$ slightly changes. This implies that this movement
of $f(x)$ inside of the $L$-type domain $D(f)$ causes a movement of
$f_e(x)=\alpha^2 (\sum_1^n k_i x_i)^2$ along the ray
$\{\lambda(\sum_1^n k_i z_i)^2: \lambda \ge 0 \}$. The collections of
the integers $\{k_i:1\le i \le n\}$ is an invariant of the $L$-type
domain $D(f)$.

\vspace{3mm}
{\bf Remark 2}. The equivalence (i)$\Leftrightarrow$(iv) of Theorem 1
is mentioned in the paper \cite{ER}. After the proof of Theorem 4.3,
the authors of \cite{ER} write:

\begin{quote}
The ideas used in this proof can be extended to cover the case in
which only a portion of the edge forms have rank 1. For such an $L$-type 
domain each rank 1 edge form can be associated with a $D$-family of 
parallel hyperplanes $G$, and the $L$-partitions $\cal S$ of lattices 
on this domain are refinement of the partition determined by $G$. In 
the other direction, any hyperplane which does not intersect the interior 
of any $L$-polytope of an $L$-partition can be associated with a rank 1 
edge form of the corresponding $L$-type domain. Such hyperplanes are 
members of a $D$-family of parallel hyperplanes.
\end{quote}

Note that here an $L$-partition and an $L$-polytope mean a Delaunay 
partition and a Delaunay polytope, and a $D$-family is a family of
parallel laminae.

In fact, the authors asserts that the ideas used in the proof of
Theorem 4.3 of \cite{ER} can be extended to the proof of the 
equivalence (i)$\Leftrightarrow$(iv).
But it seems to us that the
proof given above is not complete.
 

 \end{document}